\documentclass[12pt,a4paper]{article}

\usepackage{amsmath}
\usepackage{amsthm}
\usepackage{amssymb}
\usepackage{amscd}

\theoremstyle{definition}
\newtheorem{theorem}{Theorem}[section]
\newtheorem{prop}[theorem]{Proposition}
\newtheorem{lemma}[theorem]{Lemma}

\newtheorem{example}[theorem]{Example}

\newtheorem{conjecture}[theorem]{Conjecture}
\newtheorem*{acknowledgements}{Acknowledgements}

\newenvironment{demo}[1]{%
  \trivlist
  \item[\hskip\labelsep
        {\bf #1.}]
}{%
\hfill\qedsymbol
  \endtrivlist
}

\newcommand\Nat{\mathbb{N}}
\newcommand\Int{\mathbb{Z}}
\newcommand\ep{\varepsilon}
\newcommand\trans{{}^t\!}
\newcommand\vectx{\boldsymbol{x}}
\newcommand\vecty{\boldsymbol{y}}
\newcommand\vectz{\boldsymbol{z}}
\newcommand\tr{\operatorname{tr}}

\title{
$(q,t)$-deformations of multivariate hook product formulae
}

\author{
Soichi OKADA\thanks{
Graduate School of Mathematics, Nagoya University,
e-mail: okada@math.nagoya-u.ac.jp.
This work is partially supported by 
JSPS Grant-in-Aid for Scientific Research (C) 18540024.
}
}

\date{}

\begin{document}

\maketitle

\begin{abstract}
We generalize multivariate hook product formulae for $P$-partitions.
We use Macdonald symmetric functions to prove a $(q,t)$-deformation of Gansner's 
hook product formula for the generating functions of reverse (shifted) 
plane partitions.
(The unshifted case has also been proved by Adachi.)
For a $d$-complete poset, we present a conjectural $(q,t)$-deformation 
of Peterson--Proctor's hook product formula.
\end{abstract}

\section{
Introduction
} \label{sec:1}

R.~Stanley \cite{Stan} introduced the notion of $P$-partitions 
for a poset $P$, and studied univariate generating functions 
of them.
A {\it $P$-partition} is an order-reversing map from $P$ 
to the set of non-negative integers $\Nat$, i.e., a map 
$\sigma : P \to \Nat$ satisfying the condition:
$$
\text{if $x \le y$ in $P$, then $\sigma(x) \ge \sigma(y)$.}
$$
Let $\mathcal{A}(P)$ denote the set of all $P$-partitions.

Typical examples of $P$-partitions are reverse plane partitions 
and reverse shifted plane partitions.
If $\lambda$ is a partition, then its diagram
$$
D(\lambda) = \{ (i,j) \in \Int^2 : 1 \le j \le \lambda_i \}
$$
can be viewed as a poset by defining $(i,j) \ge (k,l)$ 
if $i \le k$ and $j \le l$, and the resulting poset is 
called a {\it shape}.
A $P$-partition for this poset $P = D(\lambda)$ is 
a {\it reverse plane partition} of shape $\lambda$,
which is an array of non-negative integers
$$
\begin{matrix}
 \pi_{1,1} & \pi_{1,2} & \multicolumn{3}{c}{\cdots\cdots}
 & \pi_{1,\lambda_1} \\
 \pi_{2,1} & \pi_{2,2} & \multicolumn{2}{c}{\cdots}
 & \pi_{2,\lambda_2} \\
 \vdots    & \vdots \\
 \pi_{r,1} & \pi_{r,2} & \cdots & \pi_{r,\lambda_r}
\end{matrix}
$$
satisfying
$$
\pi_{i,j} \le \pi_{i,j+1},
\quad
\pi_{i,j} \le \pi_{i+1,j}
$$
whenever both sides are defined.
If $\mu$ is a strict partition, then its shifted diagram
$$
S(\mu) = \{ (i,j) \in \Int^2 : i \le j \le \mu_i + i-1 \}
$$
is also a poset called a {\it shifted shape}, 
and a $S(\mu)$-partition is a {\it reverse shifted plane partition} 
of shifted shape $\mu$, which is an array of non-negative integers
$$
\begin{matrix}
 \sigma_{1,1} & \sigma_{1,2} & \sigma_{1,3}
 & \multicolumn{4}{c}{\cdots\cdots} & \sigma_{1,\mu_1} \\
              & \sigma_{2,2} & \sigma_{2,3}
 & \multicolumn{3}{c}{\cdots} & \sigma_{2,\mu_2+1} \\
              &              & \ddots \\
              &              &        & \sigma_{r,r} &
 \cdots & \sigma_{r,\mu_r+r-1}
\end{matrix}
$$
satisfying
$$
\sigma_{i,j} \le \sigma_{i,j+1},
\quad
\sigma_{i,j} \le \sigma_{i+1,j}
$$
whenever both sides are defined.

E.~Gansner \cite{Gan} considered multivariate (trace) generating 
functions for reverse (shifted) plane partitions.
Let $P$ be a shape or a shifted shape.
To each $P$-partition $\sigma$, we associate a monomial defined by
$$
\vectz^{\tr(\sigma)} = \prod_{(i,j) \in P} z_{j-i}^{\sigma(i,j)},
$$
where $z_k$ ($k \in \Int$) are indeterminates.
This weights $\sigma$ by the sums of its diagonals.
To state Gansner's hook product formulae, we introduce 
the notion of hooks for shapes and shifted shapes.
For a partition $\lambda$ and a cell $(i,j) \in D(\lambda)$, 
the {\it hook} at $(i,j)$ in $D(\lambda)$ is defined by
\begin{align}
H_{D(\lambda)}(i,j)
 &=
\{ (i,j) \} \cup \{ (i,l) \in D(\lambda) : l > j \}
 \notag
\\
 &\quad\quad
 \cup \{ (k,j) \in D(\lambda) : k > i \}.
\label{eq:hook}
\end{align}
For a strict partition $\mu$ and a cell $(i,j) \in S(\mu)$, 
the {\it shifted hook} at $(i,j)$ in $S(\mu)$ is given by
\begin{align}
H_{S(\mu)}(i,j)
 &=
 \{ (i,j) \}
 \cup
 \{ (i,l) \in S(\mu) : l > j \}
\notag
\\
 &\quad\quad
 \cup
 \{ (k,j) \in S(\mu) : k > i \}
\notag
\\
 &\quad\quad
 \cup
 \{ (j+1,l) \in S(\mu) : l > j \}.
\label{eq:shifted-hook}
\end{align}
For a finite subset $H \subset \Int^2$, we write
\begin{equation}
\vectz[H] = \prod_{(i,j) \in H} z_{j-i}.
\label{eq:monomial}
\end{equation}
Gansner \cite{Gan} used the Hillman--Grassl correspondence 
to prove the following hook product formulae.
(See also \cite{Sag}.)

\begin{theorem} \label{thm:gansner}
\cite[Theorems 5.1 and 7.1]{Gan}
Let $P$ be a shape $D(\lambda)$ or a shifted shape $S(\mu)$.
Then the multivariate generating function for $\mathcal{A}(P)$ is given by
\begin{equation}
\sum_{\sigma \in \mathcal{A}(P)} \vectz^{\tr(\sigma)}
 =
\prod_{v \in P}
 \frac{ 1 }
      { 1 - \vectz[H_P(v)] }.
\label{eq:gansner}
\end{equation}
\end{theorem}

The first aim of this paper is to prove a $(q,t)$ deformation 
of Gansner's hook product formulae.
Let $q$ and $t$ be indeterminates and put
\begin{equation}
f_{q,t}(n;m)
 =
\prod_{i=0}^{n-1}
 \frac{ 1 - q^i t^{m+1} }
      { 1 - q^{i+1} t^m }.
\label{eq:f}
\end{equation}
for non-negative integers $n$ and $m$.
And we use the notation
\begin{equation}
F(x ; q,t) = \frac{ (tx;q)_\infty }{ (x ; q)_\infty },
\label{eq:F}
\end{equation}
where $(a;q)_\infty = \prod_{i \ge 0} (1 - a q^i)$.
If we take $q=t$, then $f_{q,q}(n;m) = 1$ and $F(x;q,q) = 1/(1-x)$.
Our main theorem is the following:

\begin{theorem} \label{thm:main}
\begin{enumerate}
\item[(a)]
Let $\lambda$ be a partition.
We define a weight $W_{D(\lambda)}(\sigma;q,t)$ of a reverse plane partition
 $\pi \in \mathcal{A}(D(\lambda))$ by putting
\begin{multline}
W_{D(\lambda)}(\pi ; q,t)
\\
=
\prod_{(i,j) \in D(\lambda)}
\prod_{m \ge 0}
\frac{
 f_{q,t}( \pi_{i,j} - \pi_{i-m,j-m-1} ; m )
 f_{q,t}( \pi_{i,j} - \pi_{i-m-1,j-m} ; m )
}
{
 f_{q,t}( \pi_{i,j} - \pi_{i-m,j-m} ; m )
 f_{q,t}( \pi_{i,j} - \pi_{i-m-1,j-m-1} ; m )
},
\label{eq:weight1}
\end{multline}
where we use the convention that $\pi_{k,l} = 0$ if $k < 0$ or $l < 0$.
Then we have
\begin{equation}
\sum_{\pi \in \mathcal{A}(D(\lambda))} W_{D(\lambda)}(\pi ; q,t) \vectz^{\tr(\pi)}
 =
\prod_{v \in D(\lambda)} F( \vectz[ H_{D(\lambda)}(v) ] ; q,t ).
\label{eq:main1}
\end{equation}
\item[(b)]
Let $\mu$ be a strict partition.
We define a weight $W_{S(\mu)}(\sigma;q,t)$ of a reverse shifted plane partition
 $\sigma \in \mathcal{A}(S(\mu))$ by putting
\begin{align}
&
W_{S(\mu)}(\sigma ; q,t)
\notag\\
&=
\prod_{\substack{(i,j) \in S(\mu) \\ i < j}}
\prod_{m \ge 0}
\frac{
 f_{q,t}( \sigma_{i,j} - \sigma_{i-m,j-m-1} ; m )
 f_{q,t}( \sigma_{i,j} - \sigma_{i-m-1,j-m} ; m )
}
{
 f_{q,t}( \sigma_{i,j} - \sigma_{i-m,j-m} ; m )
 f_{q,t}( \sigma_{i,j} - \sigma_{i-m-1,j-m-1},m )
}
\notag\\
&\quad\times
\prod_{(i,i) \in S(\mu)}
\prod_{m \ge 0}
\frac{
 f_{q,t}( \sigma_{i,i} - \sigma_{i-2m-1,i-2m} ; 2m )
 f_{q,t}( \sigma_{i,i} - \sigma_{i-2m-2,i-2m-1} ; 2m+1 )
}
{
 f_{q,t}( \sigma_{i,i} - \sigma_{i-2m,i-2m} ; 2m )
 f_{q,t}( \sigma_{i,i} - \sigma_{i-2m-2,i-2m-2} ; 2m+1 )
},
\label{eq:weight2}
\end{align}
where $\sigma_{k,l} = 0$ if $k < 0$.
Then we have
\begin{equation}
\sum_{\sigma \in \mathcal{A}(S(\mu))} W_{S(\mu)}(\sigma ; q,t) \vectz^{\tr(\sigma)}
 =
\prod_{v \in S(\mu)} F( \vectz[ H_{S(\mu)}(v) ] ; q, t ).
\label{eq:main2}
\end{equation}
\end{enumerate}
\end{theorem}

Note that there are only finitely many terms different from $1$ 
in the products in (\ref{eq:weight1}) and (\ref{eq:weight2}).
If we put $q=t$, then $f_{q,q}(n;m) = 1$ and all weights 
$W_P(\sigma;q,q)$ are equal to $1$, so Theorem~{\ref{thm:main}} reduces to 
Gansner's hook product formula (Theorem~{\ref{thm:gansner}}).

S.~Adachi \cite{Ada} proves the formula (\ref{eq:main1}) by generalizing 
the argument of \cite{Vu}.
We give a similar but more transparent proof.
If $\lambda$ is a rectangular partition $(r^c)$, then reverse 
plane partitions of shape $(r^c)$ can be viewed as plane partitions 
by rotating $180^\circ$, and the formula (\ref{eq:main1}) gives Vuleti\'c's 
generalization of MacMahon's formula \cite[Theorem~A]{Vu}.
O.~Foda et al. \cite{FWZ, FW} uses fermion calculus to re-derive 
the Schur ($q=t$) and Hall--Littlewood ($q=0$) case 
of Vuleti\'c's generalization.

Gansner's hook product formulae are generalized to other posets 
than shapes and shifted shapes.
R.~Proctor \cite{Proc1, Proc2} introduced a wide class of posets, 
called {\it $d$-complete posets}, 
enjoying ``hook-length property'' and ``jeu-de-taquin property,''
as a generalization of shapes and shifted shapes.
And he announced \cite{Proc3} that, in collaboration with D.~Peterson, 
he obtained the hook product formula for $d$-complete posets, 
but their formulation and proof are still unpublished.
K.~Nakada \cite{Nak} gives a purely algebraic proof to a hook 
product formula equivalent to Peterson--Proctor's formula.
The second aim of this paper is to present a conjectural $(q,t)$-deformation 
of Peterson--Proctor's multivariate hook formula.

This paper is organized as follows.
In Section~\ref{sec:2}, we use the theory of Macdonald symmetric functions 
to give a generating function for reverse shifted plane partitions 
with prescribed shape and profile.
Section~\ref{sec:3} is devoted to the proof of our main theorem 
by using the result in Section~\ref{sec:2}.
In Section~\ref{sec:4}, we consider $d$-complete posets and give a conjectural 
$(q,t)$-deformation of Peterson--Proctor's hook formula.

\section{
Reverse shifted plane partitions of given shape and profile
} \label{sec:2}

In this section, we give a generating function for reverse 
shifted plane partitions of given shape and profile, 
from which our main theorem (Theorem~\ref{thm:main}) follows.
Our approach is based on an observation made in  
A.~Okounkov et al. \cite{OR, ORV}.

Let $\mu$ be a strict partition of length $r$ 
and $\sigma$ a reverse shifted plane partition of shape $\mu$.
We put
$$
\sigma[0]
 =
(\sigma_{r,r}, \sigma_{r-1,r-1}, \cdots, \sigma_{2,2}, \sigma_{1,1})
$$
and call it the {\it profile} of $\sigma$.
And we associate to $\sigma$ a weight $V_{S(\mu)}(\sigma ; q,t)$ 
defined by
\begin{align}
&
V_{S(\mu)}(\sigma ; q,t)
\notag
\\
&\quad=
\prod_{\substack{(i,j) \in S(\mu) \\ i < j}}
\prod_{m \ge 0}
\frac{
 f_{q,t}( \sigma_{i,j} - \sigma_{i-m,j-m-1} ; m )
 f_{q,t}( \sigma_{i,j} - \sigma_{i-m-1,j-m} ; m )
}
{
 f_{q,t}( \sigma_{i,j} - \sigma_{i-m,j-m} ; m )
 f_{q,t}( \sigma_{i,j} - \sigma_{i-m-1,j-m-1},m )
}
\notag
\\
&\quad\quad
\times
\prod_{(i,i) \in S(\mu)}
\prod_{m \ge 0}
\frac{
 f_{q,t}( \sigma_{i,i} - \sigma_{i-m-1,i-m} ; m )
}
{
 f_{q,t}( \sigma_{i,i} - \sigma_{i-m,i-m} ; m )
}.
\label{eq:weight3}
\end{align}
Here we use the convention that $\sigma_{k,l} = 0$ if $k \le 0$ or $k > l$,
 so only finitely many terms in the product over $m$ are 
different from $1$.
For a partition $\tau$, we denote by $\mathcal{A}(S(\mu), \tau)$ 
the set of reverse shifted plane partitions of shape $\mu$ 
and profile $\tau$.

\begin{theorem} \label{thm:refined}
Given strict partition $\mu$ of length $r$, let $N$ be such that $N \ge \mu_1$.
Let $\mu^c$ be the complement of $\mu$ in $[N] = \{ 1, 2, \cdots, N \}$, 
i.e.,
 $\{ \mu_1, \cdots, \mu_r \} \sqcup \{ \mu^c_1, \mu^c_2,\cdots, \mu^c_{N-r} \} = [N]$.
Let $\tau$ be a partition of length $\le r$.
Then we have
$$
\sum_{\sigma \in \mathcal{A}(S(\mu) , \tau)}
 V_{S(\mu)}(\sigma ; q,t) \vectz^{\tr(\sigma)}
=
\prod_{\mu^c_k < \mu_l}
 F( \tilde{z}_{\mu^c_k}^{-1} \tilde{z}_{\mu_l} ; q,t )
\cdot
Q_\tau( \tilde{z}_{\mu_1}, \cdots, \tilde{z}_{\mu_r}; q,t ),
$$
where the product is taken over all pairs $(k,l)$ satisfying $\mu^c_k < \mu_l$,
and
$$
\tilde{z}_0 = 1,
\quad
\tilde{z}_k = z_0 z_1 \cdots z_{k-1}
\quad(k \ge 1).
$$
And $Q_\tau(x;q,t)$ is the Macdonald symmetric function.
\end{theorem}

Here we recall the definition of Macdonald symmetric functions.
(See \cite[Chap.~VI]{Mac} for details.)
Let $\Lambda$ be the ring of symmetric functions in $x = (x_1, x_2, \cdots)$
 with coefficients in a field of characteristic $0$, which contains $q$, $t$ 
and formal power series in $z_k$'s.
Define a bilinear form on $\Lambda$ by
$$
\langle p_\lambda, p_\mu \rangle
 =
\delta_{\lambda, \mu}
z_\lambda
\prod_{i=1}^{l(\lambda)} \frac{1-q^{\lambda_i}}{1-t^{\lambda_i}},
$$
where $p_\lambda$ is the power-sum symmetric function.
%See (2.1) and (2.2)
Then \cite[Chap.~VI, (4.7)]{Mac} 
there exists a unique family of symmetric functions
 $P_\lambda(x;q.t) \in \Lambda$, indexed by partitions, satisfying the following two conditions:
\begin{enumerate}
\item[(i)]
$P_\lambda$ is a linear combination of monomial symmetric functions of 
the form
$$
P_\lambda
 = 
m_\lambda
 +
\sum_{\mu < \lambda} u_{\lambda, \mu} m_\mu,
$$
where $<$ denotes the dominance order on partitions.
\item[(ii)]
If $\lambda \neq \mu$, then $\langle P_\lambda, P_\mu \rangle = 0$.
\end{enumerate}
% See (4.7)
Then the $P_\lambda$'s form a basis of $\Lambda$.
The Macdonald $Q$ functions $Q_\lambda(x;q,t)$ are the dual basis of $P_\lambda$, i.e.,
$$
\langle P_\lambda, Q_\mu \rangle
 =
\delta_{\lambda,\mu}.
$$

We use operator calculus on the ring of symmetric functions 
$\Lambda$ to show Theorem~\ref{thm:refined}.
We denote by $g_n$ the Macdonald symmetric function $Q_{(n)}$ 
corresponding to a one-row partition $(n)$.
Let $g_n^+$ and $g_n^-$ be the multiplication and skewing operators 
on $\Lambda$ associated to $g_n$ respectively.
They satisfy
$$
g_n^+(h) = g_n h,
\quad
\langle g_n^-(h), f \rangle = \langle h, g_n f \rangle
$$
for any $f$, $h \in \Lambda$.
We consider the generating functions
$$
G^+(u) = \sum_{n=0}^\infty g_n^+ u^n,
\quad
G^-(u) = \sum_{n=0}^\infty g_n^- u^n.
$$
Also we introduce the degree operator $D(y)$ defined by
$$
D(y) P_\lambda = y^{|\lambda|} P_\lambda.
$$

For a strict partition $\mu$ and a partition $\tau$, we put
$$
R_{S(\mu), \tau}(\vectz ; q,t)
 =
\sum_{\sigma \in \mathcal{A}(S(\mu) , \tau)}
 V_{S(\mu)}(\sigma ; q,t) \vectz^{\tr(\sigma)}.
$$
Fix a positive integer $N$ satisfying $N \ge \mu_1$ and 
define a sequence $\ep = (\ep_1, \cdots, \ep_N)$ of $+$ and $-$ by putting
$$
\ep_k
 =
\begin{cases}
+ &\text{if $k$ is a part of $\mu$,} \\
- &\text{if $k$ is not a part of $\mu$.}
\end{cases}
$$
The first step of the proof of Theorem~\ref{thm:refined} is the following Lemma.

\begin{lemma} \label{lem:1}
In the ring of symmetric functions $\Lambda$, we have
\begin{multline*}
\sum_\tau
 R_{S(\mu), \tau}(\vectz ; q,t) P_\tau(x;q.t)
\\
 =
D(z_0) G^{\ep_1}(1) D(z_1) G^{\ep_2}(1) D(z_2) G^{\ep_2}(1)
 \cdots G^{\ep_{N-1}}(1) D(z_{N-1}) G^{\ep_N}(1) 1,
\end{multline*}
where $\tau$ runs over all partitions of length at most the length of $\mu$.
\end{lemma}

For example, if $\mu = (6,5,2)$ and $N=6$, then $\ep = (-,+,-,-,+,+)$ and
\begin{multline*}
\sum_\tau
 R_{S(\mu), \tau}(\vectz ; q,t) P_\tau(x;q,t)
\\
 =
D(z_0) G^-(1) D(z_1) G^+(1) D(z_2) G^-(1) D(z_3) G^-(1) D(z_4) G^+(1) D(z_5) G^+(1) 1.
\end{multline*}

\begin{demo}{Proof}
For a map $\sigma : S(\mu) \to \Nat$ (or an array of non-negative 
integers of shape $\mu$) and an integer $k$ ($0 \le k \le N$), 
we define the $k$th trace $\sigma[k]$ to be the sequence 
$(\cdots, \sigma_{2,k+2}, \sigma_{1,k+1})$ 
obtained by reading the $k$-th diagonal from SE to NW.
For example, a reverse shifted plane partition
$$
\sigma
 =
\begin{matrix}
 0 & 0 & 1 & 2 & 3 & 3 \\
   & 1 & 2 & 3 & 3 & 3 \\
   &   & 2 & 4
\end{matrix}
$$
has traces
\begin{gather*}
\sigma[0] = (2,1,0), \quad
\sigma[1] = (4,2,0), \quad
\sigma[2] = (3,1), \quad
\sigma[3] = (3,2),
\\
\sigma[4] = (3,3), \quad
\sigma[5] = (3), \quad
\sigma[6] = \emptyset.
\end{gather*}
For two partitions $\alpha$ and $\beta$, we write $\alpha \succ \beta$ 
if $\alpha_1 \ge \beta_1 \ge \alpha_2 \ge \beta_2 \ge \cdots$,
i.e, the skew diagram $\alpha/\beta$ is a horizontal strip.
Then it is clear from definition that a map $\sigma : S(\mu) \to \Nat$ 
is a shifted reverse plane partition if and only if each $\sigma[k]$ ($0 \le k \le N$)
 is a partition and
$$
\begin{cases}
\sigma[k-1] \succ \sigma[k] &\text{if $\ep_k = +$,} \\
\sigma[k-1] \prec \sigma[k] &\text{if $\ep_k = -$.}
\end{cases}
$$
In the above example, we have
$$
\sigma[0] \prec \sigma[1] \succ \sigma[2] \prec \sigma[3] 
\prec \sigma[4] \succ \sigma[5] \succ \sigma[6].
$$

A key ingredient is the Pieri rule for Macdonald symmetric functions.
For two partitions $\alpha$ and $\beta$ satisfying $\alpha \succ \beta$, 
we put
\begin{equation}
\begin{aligned}
\varphi^+_{\alpha, \beta}(q,t)
 &=
\prod_{i \le j}
 \frac{ f_{q,t}( \alpha_i - \beta_j ; j-i ) f_{q,t}( \beta_i - \alpha_{j+1} ; j-i ) }
      { f_{q,t}( \alpha_i - \alpha_j ; j-i ) f_{q,t}( \beta_i - \beta_{j+1} ; j-i ) },
\\
\varphi^-_{\beta, \alpha}(q,t)
 &=
\prod_{i \le j}
 \frac{ f_{q,t}( \alpha_i - \beta_j ; j-i ) f_{q,t}( \beta_i - \alpha_{j+1} ; j-i ) }
      { f_{q,t}( \alpha_i - \alpha_{j+1} ; j-i ) f_{q,t}( \beta_i - \beta_j ; j-i ) }.
\end{aligned}
\label{eq:phi}
\end{equation}
Then the Pieri rule \cite[Chap.~VI, (6,24)]{Mac} can be stated as follows:
\begin{equation}
\begin{aligned}
G^+(u) P_\beta
 &= 
\sum_{\alpha \succ \beta}
 \varphi^+_{\alpha, \beta}(q,t) u^{|\alpha| - |\beta|} P_\alpha,
\\
G^-(u) P_\alpha
 &= 
\sum_{\beta \prec \alpha}
 \varphi^-_{\beta, \alpha}(q,t) u^{|\alpha| - |\beta|} P_\beta,
\end{aligned}
\label{eq:pieri}
\end{equation}
where $\alpha$ in the first summation (resp. $\beta$ in the second summation) 
runs over all partition satisfying $\alpha \succ \beta$ 
(resp. $\beta \prec \alpha$).
(In \cite[Chap.~VI, (6.24)]{Mac}, the coefficients $\varphi^+_{\alpha, \beta}
 = \varphi_{\alpha, \beta}$ and $\varphi^-_{\beta, \alpha} = \psi_{\alpha, \beta}$ 
are given in terms of arm and leg lengths, but it is not hard to rewrite them 
in the form (\ref{eq:phi}) by using $f_{q,t}(m,n)$ defined by (\ref{eq:f}).
See also \cite[Chap.~VI, 6, Ex.~2]{Mac}.)
And the weight function $V_{S(\mu)}(\sigma ; q,t)$ is expressed in terms of 
$\varphi^\pm_{\alpha, \beta}(q,t)$ as
$$
V_{S(\mu)}(\sigma ;q,t)
 =
\prod_{k=1}^N
\varphi^{\ep_k}_{\sigma[k-1], \sigma[k]}(q,t).
$$

Now the claim of Lemma easily follows by induction on $\mu_1$.
\end{demo}

The second step is to rewrite the composite operators on the right-hand side of Lemma~\ref{lem:1}
 by using some commutation relations.

\begin{lemma} \label{lem:2}
Let $\mu^c$ be the strict partition formed by the complement 
of $\mu$ in $[N]$,
 i.e.,
$$
\{ \mu_1, \cdots, \mu_r \} \sqcup \{ \mu^c_1, \cdots, \mu^c_{N-r} \}
 = [N].
$$
Then we have
\begin{multline*}
D(z_0) G^{\ep_1}(1) D(z_1) G^{\ep_2}(1) D(z_2) G^{\ep_2}(1)
 \cdots G^{\ep_{N-1}}(1) D(z_{N-1}) G^{\ep_N}(1)
\\
=
\prod_{ \mu^c_k < \mu_l }
 F( \tilde{z}_{\mu^c_k}^{-1} \tilde{z}_{\mu_l} ; q,t )
\prod_{k=1}^r G^+(\tilde{z}_{\mu_k})
\prod_{l=1}^{N-r} G^-(\tilde{z}_{\mu^c_l})
D(\tilde{z}_N),
\end{multline*}
where the first product is taken over all pairs $(k,l)$ satisfying 
$\mu^c_k < \mu_l$.
\end{lemma}

\begin{demo}{Proof}
By the Pieri rule (\ref{eq:pieri}), we have
\begin{gather*}
D(z) \circ G^+ (u) = G^+ (z u) \circ D(z),
\\
D(z) \circ G^- (u) = G^- (z^{-1} u) \circ D(z).
\end{gather*}
Also we have
$$
D(z) \circ D(z') = D(z z').
$$
By using these relations, we move $D(z_i)$'s to the right and 
see that
\begin{multline*}
D(z_0) G^{\ep_1}(1) D(z_1) G^{\ep_2}(1) D(z_2) G^{\ep_2}(1)
 \cdots G^{\ep_{N-1}}(1) D(z_{N-1}) G^{\ep_N}(1)
\\
=
G^{\ep_1}(\tilde{z}_1^{\ep_1}) G^{\ep_2}(\tilde{z}_2^{\ep_2}) 
\cdots G^{\ep_N}(\tilde{z}_N^{\ep_N}) D(\tilde{z}_N),
\end{multline*}
where $z^+ =z$ and $z^- = z^{-1}$.

By the same argument as in \cite[Chap.~III, 5, Ex.~8]{Mac} 
for Hall--Littlewood functions, we can show that
$$
G^-(u) \circ G^+(v) = F(uv ; q,t) G^+(v) \circ G^-(u),
$$
where $F(x;q,t)$ is defined by (\ref{eq:F}).
(Details are left to the reader.)
It follows from this commutation relation that
\begin{multline*}
G^{\ep_1}(\tilde{z}_1^{\ep_1}) G^{\ep_2}(\tilde{z}_2^{\ep_2})
 \cdots G^{\ep_N}(\tilde{z}_N^{\ep_N})
\\
 =
\prod_{ \mu^c_k < \mu_l }
 F( \tilde{z}_{\mu^c_k}^{-1} \tilde{z}_{\mu_l} ; q,t )
\prod_{k=1}^r G^+(\tilde{z}_{\mu_k})
\prod_{l=1}^{N-r} G^-(\tilde{z}_{\mu^c_l}).
\end{multline*}
\end{demo}

Now we are in position to finish the proof of Theorem~{\ref{thm:refined}}.

\begin{demo}{Proof of Theorem~{\ref{thm:refined}}}
It follows from Lemmas~\ref{lem:1} and \ref{lem:2} that
$$
\sum_\tau
 R_{S(\mu), \tau}(\vectz ; q,t) P_\tau(x;q,t)
=
\prod_{ \mu^c_k < \mu_l }
 F( \tilde{z}_{\mu^c_k}^{-1} \tilde{z}_{\mu_l} ; q,t )
\prod_{k=1}^r G^+(\tilde{z}_{\mu_k})
\prod_{l=1}^{N-r} G^-(\tilde{z}_{\mu^c_l})
D(\tilde{z}_N) 1.
$$
By definition, we have $D(z) 1 = 1$ and $G^-(u) 1 = 1$, so we see that
$$
\sum_\tau
 R_{S(\mu), \tau}(\vectz ; q,t) P_\tau(x;q,t)
 =
\prod_{ \mu^c_k < \mu_l }
 F( \tilde{z}_{\mu^c_k}^{-1} \tilde{z}_{\mu_l} ; q,t )
\prod_{k=1}^r G^+(\tilde{z}_{\mu_k}) 1.
$$
Since the generating function of $g_n(x;q,t)$, where $x = (x_1, x_2, \cdots)$, is given by 
(see \cite[Chap.~VI, (2.8)]{Mac})
$$
\sum_{n=0}^\infty g_n(x;q,t) u^n = \prod_i F( x_i u ; q,t),
$$
we have
$$
\sum_\tau
 R_{S(\mu), \tau}(\vectz ; q,t) P_\tau (x;q,t)
=
\prod_{ \mu^c_k < \mu_l }
 F( \tilde{z}_{\mu^c_k}^{-1} \tilde{z}_{\mu_l} ; q,t )
\prod_{k=1}^r \prod_i F( x_i \tilde{z}_{\mu_k} ; q,t ).
$$
It follows from the Cauchy identity \cite[Chap.~VI, (4.13)]{Mac} that
$$
\prod_{k=1}^r \prod_i F( x_i \tilde{z}_{\mu_k} ; q,t )
 =
\sum_\tau
 Q_\tau (\tilde{z}_{\mu_1}, \cdots, \tilde{z}_{\mu_r} ; q,t) P_\tau (x;q,t).
$$
Hence we see that
\begin{multline*}
\sum_\tau
 R_{S(\mu), \tau}(\vectz ; q,t) P_\tau (x;q,t)
\\
=
\prod_{ \mu^c_k < \mu_l }
 F( \tilde{z}_{\mu^c_k}^{-1} \tilde{z}_{\mu_l} ; q,t )
\sum_\tau
 Q_\tau (\tilde{z}_{\mu_1}, \cdots, \tilde{z}_{\mu_r}  ; q,t) P_\tau (x;q,t).
\end{multline*}
Equating the coefficients of $P_\tau(x;q,t)$ completes the proof.
\end{demo}

\section{
Proof of Theorem~{\ref{thm:main}}
} \label{sec:3}

In this section, we derive Theorem~{\ref{thm:main}} from Theorem~{\ref{thm:refined}}.
If we put
$$
b_\tau (q,t)
 =
\prod_{i \le j}
 \frac{ f_{q,t}( \tau_i - \tau_{j+1} ; j-i ) }
      { f_{q,t}( \tau_i - \tau_j ; j-i ) },
$$
then we have (see \cite[Chap.~VI, (4,12) and (6.19)]{Mac})
$$
Q_\tau(x;q,t) = b_\tau(q,t) P_\tau(x;q,t).
$$

\begin{demo}{Proof of Theorem~{\ref{thm:main}} (b)}
First we prove our $(q,t)$ deformation for shifted shapes.
By comparing (\ref{eq:weight2}) and (\ref{eq:weight3}), 
for a reverse shifted plane partition $\sigma$ of shifted shape $\mu$, 
we have
$$
W_{S(\mu)}(\sigma;q,t)
 =
\frac{ b_\tau^{\text{el}}(q,t) }{ b_\tau(q,t) }
V_{S(\mu)}(\sigma;q,t),
$$
where $\tau = \sigma[0]$ is the profile of $\sigma$ and
$$
b^{\text{el}}_\tau (q,t)
 =
\prod_{\substack{i \le j \\ \text{$j-i$ is even}}}
 \frac{ f_{q,t}( \tau_i - \tau_{j+1} ; j-i ) }
      { f_{q,t}( \tau_i - \tau_j ; j-i ) }.
$$
with the product taken over all $i$ and $j$ such that $i \le j$ 
and $j-i$ is even.
Hence it follows from Theorem~{\ref{thm:refined}} that
\begin{align*}
&
\sum_{\sigma \in \mathcal{A}(S(\mu))}
 W_{S(\mu)}(\sigma ; q,t) \vectz^{\tr(\sigma)}
\\
 &\quad=
\sum_\tau \sum_{\sigma \in \mathcal{A}(S(\mu), \tau)}
 \frac{ b_\tau^{\text{el}}(q,t) }{ b_\tau(q,t) }
 V_{S(\mu)}(\sigma;q,t) \vectz^{\tr(\sigma)}
\\
 &\quad=
\sum_\tau
 \frac{ b_\tau^{\text{el}}(q,t) }{ b_\tau(q,t) }
 \prod_{\mu^c_k < \mu_l}
  F( \tilde{z}_{\mu^c_k}^{-1} \tilde{z}_{\mu_l} ; q,t)
 Q_\tau( \tilde{z}_{\mu_1}, \cdots, \tilde{z}_{\mu_r} ; q,t ).
\end{align*}
By applying the Schur--Littlewood type formula 
\cite[Chap.~VI, 7, Ex.~4 (ii)]{Mac}
$$
\sum_\tau
 \frac{ b_\tau^{\text{el}}(q,t) }{ b_\tau(q,t) }
 Q_\tau(x ; q,t)
 =
\prod_i F(x_i ; q,t)
\prod_{i < j} F(x_i x_j ; q,t),
$$
we see that
\begin{multline*}
\sum_{\sigma \in \mathcal{A}(S(\mu))}
 W_{S(\mu)}(\sigma ; q,t) \vectz^{\tr(\sigma)}
\\
 =
 \prod_{\mu^c_k < \mu_l}
  F( \tilde{z}_{\mu^c_k}^{-1} \tilde{z}_{\mu_l} ; q,t)
 \prod_{i=1}^r F(\tilde{z}_{\mu_k} ; q,t)
 \prod_{1 \le k < l \le r} F(\tilde{z}_{\mu_k} \tilde{z}_{\mu_l} ; q,t).
\end{multline*}

On the other hand, we can compute the monomial $\vectz[H_{S(\mu)}(i,j)]$ 
associated to the shifted hook $H_{S(\mu)}(i,j)$ at $(i,j) \in S(\mu)$. 
By using the fact that the length of column $j$ ($j > r$) of the shifted diagram 
$S(\mu)$ is equal to $j - \mu^c_{N-j+1}$, we have
$$
\vectz[H_{S(\mu)}(i,j)] 
 =
\begin{cases}
\tilde{z}_{\mu_i} \tilde{z}_{\mu_{j+1}} &\text{if $i \le j < r$,}
\\
\tilde{z}_{\mu_i} &\text{if $i \le j = r$,}
\\
\tilde{z}_{\mu^c_{N-j+1}}^{-1} \tilde{z}_{\mu_i} &\text{if $i \le r < j$.}
\end{cases}
$$
Noticing that $(i,j) \in S(\mu)$
if and only if $\mu_i > \mu^c_{N-j+1}$ for $1 \le i \le r < j \le N$, 
the desired product has been obtained.
\end{demo}

\begin{demo}{Proof of Theorem~{\ref{thm:main}} (a)}
Next we show the $(q,t)$-deformation for shapes.
For a given partition $\lambda$, let $r = \# \{ i : \lambda_i \ge i \}$ 
be the number of cells on the main diagonal of $D(\lambda)$ and 
define two strict partitions $\mu = (\mu_1, \cdots, \mu_r)$ and 
$\nu = (\nu_1, \cdots, \nu_r)$ by
$$
\mu_i = \lambda_i - i + 1,
\quad
\nu_i = \trans\lambda_i - i + 1
\quad(1 \le i \le r),
$$
where $\trans \lambda$ is the conjugate partition of $\lambda$.
Also we put
\begin{gather*}
x_0 = z_0^{1/2}, \qquad x_k = z_k \quad(k \ge 1),
\\
y_0 = z_0^{1/2}, \qquad y_k = z_{-k} \quad(k \ge 1).
\end{gather*}
Then a reverse plane partition $\pi \in \mathcal{A}(D(\lambda))$ is 
obtained by gluing two reverse shifted plane partitions 
$\sigma \in \mathcal{A}(S(\mu))$ and $\rho \in \mathcal{A}(S(\nu))$ 
with the same profile $\tau = \sigma[0] = \rho[0]$, and
\begin{gather*}
\vectz^{\tr(\pi)} = \vectx^{\tr(\sigma)} \vecty^{\tr(\rho)},
\\
W_{D(\lambda)}(\pi;q,t)
 =
\frac{1}{ b_\tau(q,t) } V_{S(\mu)}(\sigma;q,t) V_{S(\nu)}(\rho;q,t).
\end{gather*}
Hence it follows from Theorem~{\ref{thm:refined}} that
\begin{align*}
&
\sum_{\pi \in \mathcal{A}(D(\lambda))}
 W_{D(\lambda)}(\pi;q,t) \vectz^{\tr(\pi)}
\\
&\quad=
\sum_\tau
 \sum_{\substack{ \sigma \in \mathcal{A}(\mu,\tau) \\ \rho \in \mathcal{A}(\nu,\tau) }}
  \frac{1}{ b_\tau(q,t) } V_{S(\mu)}(\sigma;q,t) V_{S(\nu)}(\rho;q,t)
   \vectx^{\tr(\sigma)} \vecty^{\tr(\rho)}
\\
&\quad=
\sum_\tau
 \frac{1}{ b_\tau(q,t) }
 \prod_{\mu^c_k < \mu_l} F( \tilde{x}_{\mu^c_k}^{-1} \tilde{x}_{\mu_l} ; q,t )
 Q_\tau( \tilde{x}_{\mu_1}, \cdots, \tilde{x}_{\mu_r} ; q,t )
\\
&
\phantom{\quad = \sum_\tau \frac{1}{ b_\tau(q,t) } }
\times
 \prod_{\nu^c_k < \nu_l} F( \tilde{y}_{\nu^c_k}^{-1} \tilde{y}_{\nu_l} ; q,t )
 Q_\tau( \tilde{y}_{\nu_1}, \cdots, \tilde{y}_{\nu_r} ; q,t ).
\end{align*}
Now applying the Cauchy identity \cite[Chap.~VI, (4.13)]{Mac}, 
we obtain
\begin{multline*}
\sum_{\pi \in \mathcal{A}(D(\lambda))}
 W_{D(\lambda)}(\pi;q,t) \vectz^{\tr(\pi)}
\\
=
 \prod_{\mu^c_k < \mu_l} F( \tilde{x}_{\mu^c_k}^{-1} \tilde{x}_{\mu_l} ; q,t )
 \prod_{\nu^c_k < \nu_l} F( \tilde{y}_{\nu^c_k}^{-1} \tilde{y}_{\nu_l} ; q,t )
 \prod_{i,j=1}^r F( \tilde{x}_{\mu_i} \tilde{y}_{\nu_j} ; q,t ).
\end{multline*}

On the other hand, the monomial $\vectz[H_{D(\lambda)}(i,j)]$ associated to 
the hook $H_{D(\lambda)}(i,j)$ at $(i,j) \in D(\lambda)$ is given by
$$
\vectz[H_{D(\lambda)}(i,j)]
 =
\begin{cases}
\tilde{x}_{\mu_i} \tilde{y}_{\nu_j} &\text{if $1 \le i, j \le r$,}
\\
\tilde{x}_{\mu^c_{N-j+1}}^{-1} \tilde{x}_{\mu_i} &\text{if $1 \le i \le r < j$,}
\\
\tilde{y}_{\nu^c_{N-i+1}}^{-1} \tilde{y}_{\nu_j} &\text{if $1 \le j \le r < i$.}
\end{cases}
$$
This completes the proof of Theorem~{\ref{thm:main}}.
\end{demo}

Theorem~{\ref{thm:main}} (a) can be also obtained from Theorem~{\ref{thm:refined}} by 
specializing $\tau = \emptyset$.

\section{
$(q,t)$-deformation of Peterson--Proctor's hook formula
} \label{sec:4}

In this section we give a conjectural $(q,t)$-deformation of Peterson--Proctor's 
hook formulae for $d$-complete posets.

\begin{figure}[htb]
\begin{center}
\setlength{\unitlength}{2pt}
\begin{picture}(30,90)
\put(15,5){\circle*{3}}
\put(15,15){\circle*{3}}
\put(15,35){\circle*{3}}
\put(5,45){\circle*{3}}
\put(25,45){\circle*{3}}
\put(15,55){\circle*{3}}
\put(15,75){\circle*{3}}
\put(15,85){\circle*{3}}
\put(15,5){\line(0,1){10}}
\put(15,15){\line(0,1){5}}
\multiput(15,20)(0,2){5}{\line(0,1){1}}
\put(15,30){\line(0,1){5}}
\put(15,35){\line(-1,1){10}}
\put(15,35){\line(1,1){10}}
\put(5,45){\line(1,1){10}}
\put(25,45){\line(-1,1){10}}
\put(15,55){\line(0,1){5}}
\multiput(15,60)(0,2){5}{\line(0,1){1}}
\put(15,70){\line(0,1){5}}
\put(15,75){\line(0,1){10}}
\end{picture}
\end{center}
\caption{Double-tailed diamond $d_k(1)$}
\label{fig:1}
\end{figure}
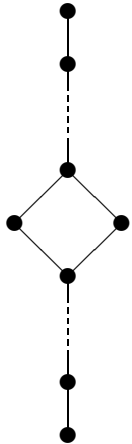

First we review the definition and some properties of $d$-complete posets.
(See \cite{Proc1, Proc2}).
For $k \ge 3$, we denote by $d_k(1)$ the poset consisting of $2k-2$ elements 
with the Hasse diagram shown in Figure~{\ref{fig:1}}.
For example, $d_3(1)$ is isomorphic to the shape $D((2,2))$,
 while $d_4(1)$ is isomorphic to the shifted shape $S((3,2,1))$.
The poset $d_k(1)$ is called the {\it double-tailed diamond} poset.

Let $P$ be a poset.
An interval $[w,v] = \{ x \in P : w \le x \le v \}$ is called 
a {\it $d_k$-interval} if it is isomorphic to $d_k(1)$.
Then $w$ and $v$ are called the {\it bottom} and {\it top} of $[w,v]$ respectively,
and the two incomparable elements of $[w,v]$ are called the {\it sides}.
If $k \ge 4$, then an interval $[w,v]$ is a {\it $d_k^-$-interval} 
if it is isomorphic 
to the poset obtained by removing the maximum element from $d_k(1)$.
A $d_3^-$-interval $[w;x,y]$ consists of three elements $x$, $y$ and $w$ 
such that $w$ is covered by $x$ and $y$.
(Precisely speaking, a $d_3^-$-interval is not an interval, 
but we use this terminology after Proctor \cite{Proc1, Proc2}.)

A poset $P$ is {\it $d$-complete} if it satisfies 
the following three conditions for every $k \ge 3$:
\begin{enumerate}
\item[(D1)]
If $I$ is a $d_k^-$-interval, then there exists an element $v$ 
such that $v$ covers the maximal elements of $I$ and 
$I \cup \{ v \}$ is a $d_k$-interval.
\item[(D2)]
If $I = [w,v]$ is a $d_k$-interval and the top $v$ covers $u$ in $P$, 
then $u \in I$.
\item[(D3)]
There are no $d_k^-$-intervals which differ only in the minimal elements.
\end{enumerate} 
It is clear that rooted trees, viewed as posets with their roots being the maximum elements, 
are $d$-complete posets.
And it can be shown that shapes and shifted shapes are $d$-complete posets.

\begin{prop} \label{prop:1} (\cite[\S 3]{Proc1})
Let $P$ be a $d$-complete poset.
Suppose that $P$ is connected, i.e., the Hasse digram of $P$ is connected.
Then we have
\begin{enumerate}
\item[(a)]
$P$ has a unique maximal element $v_0$.
\item[(b)]
For each $v \in P$, every saturated chain from $v$ 
to the maximum element $v_0$ has the same length.
Hence $P$ admits a rank function $r: P \to \Nat$ such that 
$r(x) = r(y) + 1$ if $x$ covers $y$.
\end{enumerate}
\end{prop}

Let $P$ be a poset with a unique maximal element.
The {\it top tree} $T$ of $P$ is the subgraph of 
the Hasse diagram of $P$, whose vertex set consists of all elements $x \in P$ 
such that every $y \ge x$ is covered by at most one other element.

\begin{prop} \label{prop:2} (\cite[Proposition~8.6]{Proc2})
Let $P$ be a connected $d$-complete poset and $T$ its top tree.
Let $I$ be a set of colors whose cardinality is the same as $T$.
Then a bijection $c : T \to I$ can be uniquely extended to a map 
$c : P \to I$ satisfying the following four conditions:
\begin{enumerate}
\item[(C1)]
If $x$ and $y$ are incomparable, then $c(x) \neq c(y)$.
\item[(C2)]
If $x$ covers $y$, then $c(x) \neq c(y)$.
\item[(C3)]
If an interval $[w,v]$ is a chain, then the colors $c(x)$ ($x \in [w,v]$) are 
distinct.
\item[(C4)]
If $[w,v]$ is a $d_k$-interval then $c(w) = c(v)$.
\end{enumerate}
Such a map $c : P \to I$ is called a {\it $d$-complete coloring}.
\end{prop}

\begin{example} \label{ex:1}
In the case of shapes and shifted shapes, $d$-complete colorings 
are given by contents.
\begin{enumerate}
\item[(a)]
Let $\lambda$ be a partition.
Then the top tree of the shape $D(\lambda)$ is given by
$$
T
 =
\{ (i,1) : 1 \le i \le \trans\lambda_1 \}
 \cup
\{ (1,j) : 1 \le j \le \lambda_1 \}
$$
and the content function $c : D(\lambda) \to 
\{ -\trans\lambda_1+1, \cdots, -1, 0, 1, \cdots, \lambda_1 -1 \}$ defined by
$$
c(i,j) = j-i \quad((i,j) \in D(\lambda))
$$
is a $d$-complete coloring.
\item[(b)]
If $\mu$ is a strict partition with length $\ge 2$, then the top tree 
of the shifted shape $S(\mu)$ is given by
$$
T = \{ (1,j) : 1 \le j \le \mu_1 \} \cup \{ (2,2) \},
$$
and a $d$-complete coloring $c : S(\mu) \to \{ 0, 0', 1, 2, \cdots, \mu_1-1 \}$ 
is given by
$$
c(i,j)
 =
\begin{cases}
 j-i &\text{if $i < j$,} \\
 0 &\text{if $i = j$ and $i$ is odd,} \\
 0' &\text{if $i = j$ and $i$ is even.}
\end{cases}
$$
\end{enumerate}
Figure~{\ref{fig:2}} illustrates the top trees and $d$-complete colorings of 
the shape $D((5,4,3,1))$ and the shifted shape $S((7,6,3,1))$.
The top tree consists of the nodes denoted by $\circ$.
\begin{figure}[htb]
\begin{center}
\setlength{\unitlength}{2pt}
\begin{picture}(170,85)
\put(0,15){
% shape
\put(5,15){\circle{3}}
\put(15,25){\circle{3}}
\put(25,15){\circle*{3}}
\put(25,35){\circle{3}}
\put(35,5){\circle*{3}}
\put(35,25){\circle*{3}}
\put(35,45){\circle{3}}
\put(45,15){\circle*{3}}
\put(45,35){\circle{3}}
\put(55,5){\circle*{3}}
\put(55,25){\circle{3}}
\put(65,15){\circle{3}}
\put(75,5){\circle{3}}
\put(6,16){\line(1,1){8}}
\put(16,24){\line(1,-1){8}}
\put(16,26){\line(1,1){8}}
\put(26,14){\line(1,-1){8}}
\put(26,16){\line(1,1){8}}
\put(26,34){\line(1,-1){8}}
\put(26,36){\line(1,1){8}}
\put(36,6){\line(1,1){8}}
\put(36,24){\line(1,-1){8}}
\put(36,26){\line(1,1){8}}
\put(36,44){\line(1,-1){8}}
\put(46,14){\line(1,-1){8}}
\put(46,16){\line(1,1){8}}
\put(46,34){\line(1,-1){8}}
\put(56,6){\line(1,1){8}}
\put(56,24){\line(1,-1){8}}
\put(66,14){\line(1,-1){8}}
\put(0,18){\makebox(10,5){$-3$}}
\put(10,28){\makebox(10,5){$-2$}}
\put(20,18){\makebox(10,5){$-1$}}
\put(20,38){\makebox(10,5){$-1$}}
\put(30,8){\makebox(10,5){$0$}}
\put(30,28){\makebox(10,5){$0$}}
\put(30,48){\makebox(10,5){$0$}}
\put(40,18){\makebox(10,5){$1$}}
\put(40,38){\makebox(10,5){$1$}}
\put(50,8){\makebox(10,5){$2$}}
\put(50,28){\makebox(10,5){$2$}}
\put(60,18){\makebox(10,5){$3$}}
\put(70,8){\makebox(10,5){$4$}}
}
\put(100,0){
% shifted shape
\put(5,15){\circle*{3}}
\put(5,35){\circle*{3}}
\put(5,55){\circle{3}}
\put(5,75){\circle{3}}
\put(15,25){\circle*{3}}
\put(15,45){\circle*{3}}
\put(15,65){\circle{3}}
\put(25,15){\circle*{3}}
\put(25,35){\circle*{3}}
\put(25,55){\circle{3}}
\put(35,25){\circle*{3}}
\put(35,45){\circle{3}}
\put(45,15){\circle*{3}}
\put(45,35){\circle{3}}
\put(55,5){\circle*{3}}
\put(55,25){\circle{3}}
\put(65,15){\circle{3}}
\put(6,16){\line(1,1){8}}
\put(6,34){\line(1,-1){8}}
\put(6,36){\line(1,1){8}}
\put(6,54){\line(1,-1){8}}
\put(6,56){\line(1,1){8}}
\put(6,74){\line(1,-1){8}}
\put(16,24){\line(1,-1){8}}
\put(16,26){\line(1,1){8}}
\put(16,44){\line(1,-1){8}}
\put(16,46){\line(1,1){8}}
\put(16,64){\line(1,-1){8}}
\put(26,16){\line(1,1){8}}
\put(26,34){\line(1,-1){8}}
\put(26,36){\line(1,1){8}}
\put(26,54){\line(1,-1){8}}
\put(36,24){\line(1,-1){8}}
\put(36,26){\line(1,1){8}}
\put(36,44){\line(1,-1){8}}
\put(46,14){\line(1,-1){8}}
\put(46,16){\line(1,1){8}}
\put(46,34){\line(1,-1){8}}
\put(56,6){\line(1,1){8}}
\put(56,24){\line(1,-1){8}}
\put(0,18){\makebox(10,5){$0'$}}
\put(0,38){\makebox(10,5){$0$}}
\put(0,58){\makebox(10,5){$0'$}}
\put(0,78){\makebox(10,5){$0$}}
\put(10,28){\makebox(10,5){$1$}}
\put(10,48){\makebox(10,5){$1$}}
\put(10,68){\makebox(10,5){$1$}}
\put(20,18){\makebox(10,5){$2$}}
\put(20,38){\makebox(10,5){$2$}}
\put(20,58){\makebox(10,5){$2$}}
\put(30,28){\makebox(10,5){$3$}}
\put(30,48){\makebox(10,5){$3$}}
\put(40,18){\makebox(10,5){$4$}}
\put(40,38){\makebox(10,5){$4$}}
\put(50,8){\makebox(10,5){$5$}}
\put(50,28){\makebox(10,5){$5$}}
\put(60,18){\makebox(10,5){$6$}}
}
\end{picture}
\end{center}
\caption{Top trees and $d$-complete colorings}
\label{fig:2}
\end{figure}
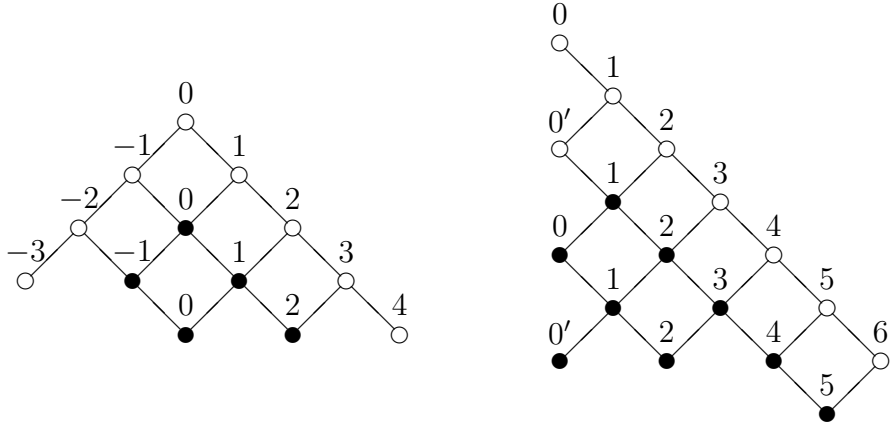
\end{example}

Let $P$ be a connected $d$-complete poset and 
$c : P \to I$ a $d$-complete coloring.
Let $z_i$ ($i \in I$) be indeterminates.
For a $P$-partition $\sigma \in \mathcal{A}(P)$, we put
$$
\vectz^\sigma = \prod_{v \in P} z_{c(v)}^{\sigma(v)}.
$$
For example, if we use the $d$-complete colorings given in Example~\ref{ex:1}, 
the monomial $\vectz^{\tr(\sigma)}$ is the same as $\vectz^\sigma$ 
for a reverse plane partition $\sigma \in \mathcal{A}(D(\lambda))$, 
while, for a reverse shifted plane partition $\sigma \in \mathcal{A}(S(\mu))$, 
$\vectz^{\tr(\sigma)}$ is obtained from $\vectz^\sigma$ by putting $z_0 = z_{0'}$.
Instead of giving a definition of hooks $H_P(v)$ for a general 
$d$-complete poset $P$, we define associated monomials $\vectz[H_P(v)]$ 
directly by induction as follows:
\begin{enumerate}
\item[(a)]
If $v$ is not the top of any $d_k$-interval, then we define
$$
\vectz[H_P(v)] = \prod_{w \le v} z_{c(w)}.
$$
\item[(b)]
If $v$ is the top of a $d_k$-interval $[w,v]$, then we define
$$
\vectz[H_P(v)]
 =
\frac{ \vectz[H_P(x)] \cdot \vectz[H_P(y)] }
     { \vectz[H_P(w)] },
$$
where $x$ and $y$ are the sides of $[w,v]$.
\end{enumerate}
It is easy to see that, for shapes, this definition of $\vectz[H_{D(\lambda)}(v)]$ 
is consistent with the definition (\ref{eq:hook}) and (\ref{eq:monomial}) 
given in Introduction.
And, for shifted shapes, the monomial $\vectz[H_{S(\mu)}(v)]$ defined above reduces
 to the monomial $\vectz[H_{S(\mu)}(v)]$ defined by (\ref{eq:shifted-hook}) 
and (\ref{eq:monomial}), if we put $z_0 = z_{0'}$.

Now we are ready to state our conjectural $(q,t)$-deformation of 
Peterson--Proctor's hook formula.
Let $P$ be a connected $d$-complete poset with the maximum element $v_0$, 
and the rank function $r : P \to \Nat$.
Let $T$ be the top tree of $P$.
Take $T$ as a set of colors and $c : P \to T$ be the $d$-complete coloring
 such that $c(v) = v$ for all $v \in T$.
Given a $P$-partition $\sigma \in \mathcal{A}(P)$, we define 
a weight $W_P(\sigma ; q,t)$ by putting
\begin{multline}
W_P(\sigma;q,t)
\\
 =
\frac{
 \displaystyle
 \prod_{\substack{x, y \in P \\ \text{s.t. } x < y, \, c(x) \sim c(y)}}
  f_{q,t} ( \sigma(x) - \sigma(y) ; d(x,y) )
 \displaystyle
 \prod_{\substack{x \in P \\ \text{s.t. } c(x) = v_0}}
  f_{q,t} ( \sigma(x); e(x,v_0) )
}
{
 \displaystyle
 \prod_{\substack{x, y \in P \\ \text{s.t. } x < y, \, c(x) = c(y)}}
 f_{q,t} ( \sigma(x) - \sigma(y) ; e(x,y) )
 f_{q,t} ( \sigma(x) - \sigma(y) ; e(x,y) - 1 )
},
\label{eq:weight}
\end{multline}
where $c(x) \sim c(y)$ means that $c(x)$ and $c(y)$ are adjacent in $T$,
and
$$
d(x,y) = \frac{ r(y) - r(x) - 1}{2},
\quad
e(x,y) = \frac{ r(y) - r(x) }{2}.
$$
Note that, if $c(x) \sim c(y)$, then $r(y) - r(x)$ is odd, and, 
if $c(x) = c(y)$, then $r(y) - r(x)$ is even, 
hence $d(x,y)$ and $e(x,y)$ are integers.

If we consider the extended poset $\hat{P} = P \sqcup \{ \hat{1} \}$, 
and its top tree $\hat{T} = T \sqcup \{ \hat{1} \}$, 
where $\hat{1}$ is the new maximum element of $\hat{P}$ and 
$\hat{1}$ is adjacent to $v_0$ in $\hat{T}$, 
then the weight $W_P(\sigma;q,t)$ can be expressed in the following form:
\begin{multline}
W_P(\sigma;q,t)
\\
 =
\frac{
 \displaystyle
 \prod_{\substack{x, y \in \hat{P} \\ \text{s.t. } x < y, \, \hat{c}(x) \sim \hat{c}(y)}}
  f_{q,t} ( \hat{\sigma}(x) - \hat{\sigma}(y) ; d(x,y) )
}
{
 \displaystyle
 \prod_{\substack{x, y \in P \\ \text{s.t. } x < y, \, c(x) = c(y)}}
 f_{q,t} ( \sigma(x) - \sigma(y) ; e(x,y) )
 f_{q,t} ( \sigma(x) - \sigma(y) ; e(x,y) - 1 )
},
\end{multline}
where $\hat{c} : \hat{P} \to \hat{T}$ and $\hat{\sigma} : \hat{P} \to \Nat$ 
are the extensions of $c$ and $\sigma$ defined by
$\hat{c}(\hat{1}) = \hat{1}$ and $\hat{\sigma}(\hat{1}) = 0$ respectively.

\begin{conjecture} \label{conj}
Let $P$ be a connected $d$-complete poset.
Using the notations defined above, we have
\begin{equation}
\sum_{\sigma \in \mathcal{A}(P)}
 W_P(\sigma ; q,t) \vectz^\sigma
 =
\prod_{v \in P}
 F( \vectz[H_P(v)] ; q,t).
\label{eq:conj}
\end{equation}
\end{conjecture}

If we put $q=t$, then $W_P(\sigma;q,q) = 1$ and $F(x;q,q) = 1/(1-x)$, 
hence Conjecture~{\ref{conj}} reduces to Peterson--Proctor's hook formula 
(see \cite{Nak}).

Our main theorem (Theorem~{\ref{thm:main}}) supports Conjecture~{\ref{conj}}.

\begin{prop} \label{prop:3}
\begin{enumerate}
\item[(a)]
Conjecture~{\ref{conj}} is true for shapes.
\item[(b)]
Conjecture~{\ref{conj}} is true for shifted shapes.
\end{enumerate}
\end{prop}

\begin{demo}{Proof}
First we check that the expressions (\ref{eq:weight1}) 
and (\ref{eq:weight2}) coincide with (\ref{eq:weight}).
Let $\sigma$ be a reverse shifted plane partition of shape $\mu$.
Then, by noting the convention that $\sigma_{k,l} = 0$ for $k < 0$, 
the expression (\ref{eq:weight2}) can be rewritten as
\begin{align*}
&
W_{S(\mu)}(\sigma;q,t)
\\
&\quad=
\prod_{\substack{(i,j) \in S(\mu) \\ i < j}}
 \prod_{m=0}^{i-2}
\frac{
 f_{q,t}( \sigma_{i,j} - \sigma_{i-m,j-m-1} ; m )
 f_{q,t}( \sigma_{i,j} - \sigma_{i-m-1,j-m} ; m )
}
{
 f_{q,t}( \sigma_{i,j} - \sigma_{i-m,j-m} ; m )
 f_{q,t}( \sigma_{i,j} - \sigma_{i-m-1,j-m-1},m )
}
\\
&\phantom{
\quad=\prod_{\substack{(i,j) \in S(\mu) \\ i < j}}
}\times
\frac{
 f_{q,t}( \sigma_{i,j} - \sigma_{1,j-i} ; i-1 )
}
{
 f_{q,t}( \sigma_{i,j} - \sigma_{1,j-i+1} ; i-1 )
}
\\
&\quad\quad\times
\prod_{\substack{(i,i) \in S(\mu) \\ \text{$i$ is odd}}}
\prod_{m=0}^{(i-3)/2}
\frac{
 f_{q,t}( \sigma_{i,i} - \sigma_{i-2m-1,i-2m} ; 2m )
 f_{q,t}( \sigma_{i,i} - \sigma_{i-2m-2,i-2m-1} ; 2m+1 )
}
{
 f_{q,t}( \sigma_{i,i} - \sigma_{i-2m,i-2m} ; 2m )
 f_{q,t}( \sigma_{i,i} - \sigma_{i-2m-2,i-2m-2} ; 2m+1 )
}
\\
&\phantom{
\quad\quad\times
\prod_{\substack{(i,i) \in S(\mu) \\ \text{$i$ is odd}}}
}\times
\frac{
 f_{q,t}( \sigma_{i,i} ; i-1 )
}
{
 f_{q,t}( \sigma_{i,i} - \sigma_{1,1} ; i-1 )
}
\\
&\quad\quad\times
\prod_{\substack{(i,i) \in S(\mu) \\ \text{$i$ is even}}}
\prod_{m=0}^{i/2-2}
\frac{
 f_{q,t}( \sigma_{i,i} - \sigma_{i-2m-1,i-2m} ; 2m )
 f_{q,t}( \sigma_{i,i} - \sigma_{i-2m-2,i-2m-1} ; 2m+1 )
}
{
 f_{q,t}( \sigma_{i,i} - \sigma_{i-2m,i-2m} ; 2m )
 f_{q,t}( \sigma_{i,i} - \sigma_{i-2m-2,i-2m-2} ; 2m+1 )
}
\\
&\phantom{
\quad\quad\times
\prod_{\substack{(i,i) \in S(\mu) \\ \text{$i$ is even}}}
}\times
\frac{
 f_{q,t}( \sigma_{i,i} - \sigma_{1,2} ; i-2 )
}
{
 f_{q,t}( \sigma_{i,i} - \sigma_{2,2} ; i-2 )
}.
\end{align*}
This expression can be transformed into (\ref{eq:weight}).
Similarly, we can show that the expression (\ref{eq:weight1}) 
is the same as (\ref{eq:weight}) for shapes.

Now (a) follows from Theorem~{\ref{thm:main}} (a).
To prove (b), we use the $d$-complete coloring given in 
Example~{\ref{ex:1}} (b) 
and modify the proof of Theorem~{\ref{thm:main}} (b), which corresponds 
to the case of $z_0 = z_{0'}$. 
(This modification is due to M.~Ishikawa \cite{I}.)
For a partition $\tau$, we denote by $o(\tau)$ the number of 
columns of odd length in the diagram of $\tau$.
Since we have
$$
\tau_1 + \tau_3 + \cdots
 = \frac{1}{2} (|\tau| + o(\tau)),
\quad
\tau_2 + \tau_4 + \cdots
 = \frac{1}{2} (|\tau| - o(\tau)),
$$
the monomial $\vectz^\sigma$ associated to $\sigma \in \mathcal{A}(S(\mu))$ 
is given by
$$
\vectz^\sigma
 =
\begin{cases}
\left( \dfrac{z_{0'}}{z_0} \right)^{(|\tau|+o(\tau))/2}
 \vectz^{\tr(\sigma)}
 &\text{if $l(\mu)$ is even,} \\
\left( \dfrac{z_{0'}}{z_0} \right)^{(|\tau|-o(\tau))/2}
 \vectz^{\tr(\sigma)}
 &\text{if $l(\mu)$ is odd,}
\end{cases}
$$
where $\tau = \sigma[0]$ is the profile of $\sigma$.
We appeal to Warnaar's generalization of the Schur--Littlewood type identity 
\cite[(1.18)]{W}
$$
\sum_\tau
 a^{o(\tau)}
 \frac{ b^{\text{el}}_\tau(q,t) }
      { b_\tau(q,t) }
 Q_\tau(x;q,t)
 =
\prod_i F(a x_i ; q.t)
\prod_{i<j} F( x_i x_j ; q,t ).
$$
Replacing $a$, $x_i$ by $a^{1/2}$, $a^{1/2} x_i^{1/2}$, or
$a^{-1/2}$, $a^{1/2}x_i^{1/2}$, respectively, we obtain
\begin{align*}
\sum_\tau
 a^{(|\tau|+o(\tau))/2}
 \frac{ b^{\text{el}}_\lambda(q,t) }
      { b_\lambda(q,t) }
 Q_\tau(x;q,t)
 &=
\prod_i F(ax_i;q,t)
\prod_{i<j} F(ax_ix_j;q,t),
\\
\sum_\tau
 a^{(|\tau|-o(\tau))/2}
 \frac{ b^{\text{el}}_\lambda(q,t) }
      { b_\lambda(q,t) }
 Q_\tau(x;q,t)
 &=
\prod_i F(x_i;q,t)
\prod_{i<j} F(a x_i x_j ; q,t).
\end{align*}
We use these identities with $a = z_{0'}/z_0$ together with Theorem~{\ref{thm:refined}} 
to derive (\ref{eq:conj}) for $P = S(\mu)$, by the argument 
similar to the proof of Theorem~{\ref{thm:main}} (b).
\end{demo}

\begin{prop} \label{prop:4}
Conjecture~{\ref{conj}} is true for rooted trees.
\end{prop}

\begin{demo}{Proof}
Let $T$ be a rooted tree with root $v_0$.
Then $\vectz[H_T(v)] = \prod_{w \le v} z_w$, and the weight function 
(\ref{eq:weight}) reduces to
$$
W_T(\sigma;q,t)
 =
\prod_{\substack{x,y \in P \\ \text{s.t. $y$ covers $x$}}}
 f_{q,t}(\sigma(x) - \sigma(y) ; 0)
\cdot
f_{q,t} (\sigma(v_0) ; 0),
$$
because $T$ itself is the top tree.

We proceed by induction on the number of vertices in the tree $T$.
Since it is obvious for $\# T = 1$, we may assume $\# T > 1$.
Let $v_1, \cdots, v_r$ be the children of $v_0$, and $T_1, \cdots, T_r$
 be the rooted subtrees of $T$ with roots $v_1, \cdots, v_r$ respectively.
Then there is a natural bijection
$$
S: \mathcal{A}(T_1) \times \cdots \times \mathcal{A}(T_r) \times \Nat
 \longrightarrow \mathcal{A}(T),
$$
which associates to $(\pi_1, \cdots, \pi_r, k) \in 
\mathcal{A}(T_1) \times \cdots \times \mathcal{A}(T_r) \times \Nat$, 
the $T$-partition $\pi \in \mathcal{A}(T)$ defined by
$$
\pi(v)
 =
\begin{cases}
 \pi_i(v)+k &\text{if $v \in T_i$}, \\
 k &\text{if $v = v_0$}.
\end{cases}
$$
Under this bijection $S$, we have
$$
W_T(\pi ; q,t)
 =
f_{q,t}(k;0) \prod_{i=1}^r W_{T_i}(\pi_i ; q,t),
\quad
\vectz^{\pi} = \vectz[T]^{k} \prod_{i=1}^r \vectz^{\pi_i}.
$$
Hence we have
\begin{multline*}
\sum_{\pi \in \mathcal{A}(T)}
 W_T(\pi ; q,t) \vectz^{\pi}
 =
\left( \sum_{k \ge 0} f_{q,t}(k;0) \vectz[T]^k \right)
\cdot
\prod_{i=1}^r 
 \left( \sum_{\pi_i \in \mathcal{A}(T_i)} W_{T_i}(\pi_i ; q,t) \vectz^{\pi_i} \right).
\end{multline*}
By using the binomial theorem
$$
\sum_{k \ge 0} f_{q,t}(k;0) \vectz[T]^k
 =
F( \vectz[T] ; q,t)
 =
F( \vectz[H_T(v_0)] ; q,t),
$$
and the induction hypothesis
$$
\sum_{\pi_i \in \mathcal{A}(T_i)} W_{T_i}(\pi_i ; q,t) \vectz^{\pi_i}
 =
\prod_{v \in T_i} F( \vectz[H_{T_i}(v)] ; q,t ),
% =
%\prod_{v \in T_i} F( \vectz[H_T(v)] ; q,t ),
$$
we can complete the proof.
\end{demo}

\begin{acknowledgements}
The author is grateful to Professor Masao Ishikawa for providing a proof of 
Proposition~\ref{prop:3} (b).
\end{acknowledgements}


\begin{thebibliography}{99}

\bibitem{Ada}
S.~Adachi,
Reverse plane partitions and growth diagrams,
preprint.

\bibitem{FW}
O.~Foda and M.~Wheeler,
Hall-Littlewood plane partitions and KP,
Int. Math. Res. Not. IMRN {\bf 2009} (2009) 2597--2619.

\bibitem{FWZ}
O.~Foda, M.~Wheeler, and M.~Zuparik,
On free fermions and plane partitions,
J. Algebra {\bf 321} (2009), 3249--3273.

\bibitem{Gan}
E.~R.~Gansner, 
The Hillman-Grassl correspondence and the enumeration of reverse plane partitions. 
J. Combin. Theory Ser. A {\bf 30} (1981), 71--89. 

\bibitem{I}
M.~Ishikawa,
private communication dated Sep. 17, 2009.

\bibitem{Mac}
I.~G.~Macdonald,
Symmetric Functions and Hall Polynomials (2nd ed.), 
Oxford Univ. Press, 1995.

\bibitem{Nak}
K.~Nakada,
$q$-Hook formula for a generalized Young diagram, 
preprint.

\bibitem{OR}
A.~Okounkov and N.~Reshetikhin,
Correlation function of Schur process with application 
to local geometry of a random 3-dimensional Young diagram,
J. Amer. Math. Soc. {\bf 16} (2003), 581--603.

\bibitem{ORV}
A.~Okounkov, N.~Reshetikhin, and C.~Vafa,
Quantum Calabi-Yau and classical crystals,
in ``The unity of mathematics'', Progr. Math., 244, 
Birkh\"auser Boston, 2006, pp. 597--618,

\bibitem{Proc1}
R.~A.~Proctor, 
Dynkin diagram classification of $\lambda$-minuscule Bruhat lattices 
and $d$-complete posets,
J. Algebraic Combin. {\bf 9} (1999), 61--94.

\bibitem{Proc2}
R.~A.~Proctor, 
Minuscule elements of Weyl groups, the number game, and $d$-complete posets,
J. Algebra {\bf 213} (1999), 272--303.

\bibitem{Proc3}
R.~A.~Proctor, 
Informal description of the hook length property for posets,
{\tt http://www.math.unc.edu/Faculty/rap/Hook.html}.

\bibitem{Sag}
B.~E.~Sagan,
Combinatorial proofs of hook generating functions for skew plane partitions,
Theoret. Comput. Sci. {\bf 117} (1993), 273--287.

\bibitem{Stan}
R.~P.~Stanley,
Ordered structures and partitions,
Mem. Amer. Math. Soc. {\bf 119} (1972).

\bibitem{Vu}
M.~Vuleti\'c,
A generalization of MacMahon's formula,
Trans. Amer. Math. Soc. {\bf 361} (2009), 2789--2804.

\bibitem{W}
S.~O.~Warnaar,
Rogers-Szeg\"o polynomials and Hall--Littlewood symmetric functions,
J. Algebra {\bf 303} (2006), 810--830.

\end{thebibliography}
\end{document}